\documentclass[12pt]{article}
\usepackage{amssymb,amsmath}
\usepackage[symbol]{footmisc}
\title{Two analogs of Pleijel's inequality}
\author{Sergey Y. Sadov}
\date{}
\newtheorem{thm}{Theorem}
\newcommand{\St}{Stieltjes transform }
\newcommand{\RM}{Riesz mean}
\newcommand{\lam}{\lambda}
\newcommand{\Nl}{N(\lam)}
\newcommand{\RR}{\mathbb{R}}
\newcommand{\CC}{\mathbb{C}}
\newcommand{\iu}{\int_0^u}
\newcommand{\iuu}{\int_{-u}^u}
\newcommand{\iR}{\int_0^{\infty}}
\newcommand{\iC}{\int_{\Gamma}}
\newcommand{\ip}{\frac{1}{2\pi i}}

\newcommand{\al}{\alpha}

\newcommand{\ze}{\zeta}
\newcommand{\eps}{\varepsilon}
\newcommand{\zbar}{\overline{\zeta}}
\newcommand{\lz}{\lam_0}
\newcommand{\zz}{\ze_0}
\newcommand{\ez}{\eta_0}
\newcommand{\tap}[1]{\tau^{\alpha#1}}
\newcommand{\tsq}{1+\tau^2}
\newcommand{\tsqq}{(\tsq)^2}
\newcommand{\f}{\frac}
\newcommand{\intdR}{S(\ze)\left(1-\f{\ze}{\lz}\right)^\al
    \,d\ze}
\newcommand{\Rzz}{R_{\al}(\zz)}
\newcommand{\sqcc}{\sqrt{|c_1|^2+c_2^2}}
\newcommand{\elz}{\left(\f{\ez}{\lz}\right)^\al}
\newcommand{\epa}{e^{i\pi\al/2}}
\newcommand{\enpa}{e^{-i\pi\al/2}}
\newcommand{\sgn}{\mathrm{sgn}\,}
\renewcommand{\Re}{\mathrm{Re}\,}
\renewcommand{\Im}{\mathrm{Im}\,}
\newcommand{\dst}{\displaystyle}
\newcommand{\beq}[1]{\begin{equation}\label{#1}}
\newcommand{\eeq}{\end{equation}}
\newcommand{\ba}{\begin{array}}
\newcommand{\ea}{\end{array}}

\DefineFNsymbols{FNstars}[text]{{} {*} {*} {**}}
\setfnsymbol{FNstars}

\begin{document}
\maketitle

\section{Formulation of results}

Let $\Nl$ be a nondecreasing function defined on $\RR_+=(0,+\infty)$
such that $\Nl=0$ for small $\lam$ and
%----------------------
\footnote{%\footnotesize\rm
Translated by the author from the Russian original: 
{\it Funkcionalnyj Analiz. Spektral'naja teorija.}
Mezhvuzovskij sbornik nauchnyh trudov.
Ul'anovsk, 1987, pp.\ 156--164.
(MR 92j:26013)}
%----------------------
\beq{eq1}
\iR\lam^{-1}d\Nl <\infty.
\eeq
The {\it \St} of $\Nl$ is defined as
\footnote{There exists an alternative convention 
according to which the \St of $f(t)$ is defined as
$\int_0^{\infty} (\ze+t)^{-1} f(t)\,dt.$}
\beq{eq2}
S(\ze)=\iR (\lam-\ze)^{-1}\,d\Nl,\qquad \ze\not\in\RR_+.
\eeq
Fix a point $\zz=\lz+i\ez\;$ in the first quadrant of the complex
plane $\CC$. Denote by $\Gamma$ a contour that connects the point
$\zz$ to $\zbar_0=\lz-i\ez$ and does not cross the integration path $\RR_+$ of (\ref{eq2}).

\AA{}ke~Pleijel in \cite{Pleijel63} obtained the inequality
\beq{eq3}
\left|N(\lz)| -\ip \iC S(\ze)\,d\ze\right|
\leq \ez\, \sqrt{1+\pi^{-2}} \;|S(\zz)|
\eeq
and used it to give a short proof of Malliavin's 
\cite{Malliavin} Tauberian theorem.
The inequality (\ref{eq3}) found applications in spectral
theory of differential and pseudo-differential operators
(see e.g.\ \cite{Agmon, PhamTheLai}).
In the present paper two generalizations of this 
inequality are derived.

\smallskip
For $\alpha>0$, the {\it \RM\ of order $\alpha$ of $\Nl$} is \beq{eq4}
N^{(\alpha)}(\lam)=\int_0^{\lam}\,\left(1-\f{x}{\lam}\right)^\al
\,dN(x),\qquad \lam>0.
\eeq
Given a power asymptotics of the \St of $\Nl$ along a certain parabola-like curve in $\CC$ that avoids $\RR_+$, the asymptotics of the \RM{}s as $\lam\to+\infty$ can be recovered using the following theorem.

\begin{thm}
\label{thm1}
Let the function $\Nl$ be constant in a neigbour\-hood of $\lz$.
Then for any $\alpha>0$
\beq{eq5}
\left|N^{(\alpha)}(\lz)| -\ip \iC \intdR\right|
\leq \f{1}{\alpha\pi}\elz\ez \,|S(\zz)|.
\eeq
For $\alpha<1$ the factor $(\alpha\pi)^{-1}$ in the right-hand
side may be replaced by $\sqrt{\pi^{-2}+1/4}$. 
\footnote{(Added in translation.)\@ 
A continuous dependence of the constant on $\al$ can be 
achieved by replacing the term $1/4$ with $(1-\al^{1+\eps})/4$, where $\eps\geq\eps_0\approx 1/16$ (found numerically).   
}
\end{thm}

\noindent
Henceforth the branch
$z^\al=\exp(\alpha\ln z)$ with $-\pi<\mathrm{Im}\,\ln z \leq\pi$ is assumed.

\smallskip
If instead of (\ref{eq1}) a weaker%
\footnote{Inequality (\ref{eq6}) is weaker than (\ref{eq1}) since we assume that $d\Nl=0$ near $\lam=0$.}  
condition with some integer $q>1$
\beq{eq6}
\iR\lam^{-q}d\Nl <\infty
\eeq
holds, then
the leading term of the asymptotics of $\Nl$ can be recovered
by means of the next theorem from the behaviour of its {\it generalized \St}
\beq{eq7}
S_q(\ze)=\iR (\lam-\ze)^{-q}\,d\Nl,\qquad \ze\not\in\RR_+.
\eeq
\begin{thm}
\label{thm2}
Let the function $\Nl$ satisfy {\rm(\ref{eq6})} and
be constant in a neigbourhood of $\lz$.
There exist constants
$C_0, C_1,\dots,$ $C_{q-2}$ (which depend only on $q$) such that
\beq{eq8}
\ba{l}
\dst
\left|N(\lz)-\ip \iC S_q(\ze)(\ze-\lz)^{q-1}\,d\ze\right|
\\[3ex]
\hspace*{5em}\dst
\leq\sum_{m=0}^{q-2}\,C_m\ez^{q-1-m}\,\left|\iC S_q(\ze)(\lz-\ze)^m\,
d\ze\right|.
\ea
\eeq
\end{thm}

%---------------------------------------------------------
\section{Proof of Theorem~\ref{thm1}}

{\bf a}. The left-hand side of (\ref{eq5}) vanishes if one uses a closed contour of integration consisting of $\Gamma$ and the segment $[\ze,\bar{\ze_0}]$.
Indeed, since $\Nl$ is assumed constant in the vicinity of
$\lz$, one may change the order of integration:
$$
\ip\oint\intdR=\iR\,d\Nl\,\ip\oint\left(\f{\lz-\ze}{\lz}\right)^\al
\,\f{d\ze}{\lam-\ze}.
$$
The inner integral in the r.h.s.\ equals $2\pi i(1-\lam/\lz)^\al$
when $\lam<\lz$, and $0$ when $\lam>\lz$. Thus the r.h.s.\
equals $N^{(\alpha)}(\lz)$. In order to prove (\ref{eq5}) we have to evaluate
\beq{eq9}
\Rzz=\ip\int\limits_{\zbar_0}^{\ze_0}\intdR
\eeq
(where the integration path is the vertical segment).

\smallskip
Set $\,\ze=\lz+i|\lam-\lz|\tau$, $s(\lam)=\sgn(\lam-\lz)$,
$u(\lam)=\ez|\lam-\lz|^{-1}$. Changing the order of integration in
(\ref{eq9}), we get
\beq{eq10}
\Rzz=\f{(-i)^\al}{2\pi}\elz\,
\iR\f{d\Nl}{u^\al(\lam)}\,\int\limits_{-u(\lam)}^{u(\lam)}\f{%
\tau^\al\,\big(s(\lam)+i\tau\big)}{\tsq}\,d\tau.
\eeq
We will find constants $c_1\in\CC$ and $c_2>0$ so that for any  $u>0$ and $s=\pm1$ the inequality
\beq{eq11}
\left|\f{1}{u^\al}\,\iuu\f{\tau^\al\,(s+i\tau)}{\tsq}\,d\tau \,-\,c_1\,\f{su}{1+u^2}
\right|\;\leq\;
c_2\,\f{u^2}{1+u^2}
\eeq
will hold.
Having (\ref{eq11}) and the identities
$$
\iR\f{u(\lam) s(\lam)}{1+u^2(\lam)}\,d\Nl\,=\,\ez\,\Re S(\zz)
$$
and
$$
\iR\f{u^2(\lam)}{1+u^2(\lam)}\,d\Nl\,=\,\ez\,\Im S(\zz)
$$
one can estimate the r.h.s.\ in (\ref{eq10}) by means of the Schwarz inequality:
\beq{eq12}
\ba{lcl}\dst
 |\Rzz|&\leq&\dst \elz\,\f{\ez}{2\pi}\;\sqcc
\;\;\sqrt{\Re^2 S(\zz)+\Im^2 S(\zz)}
\\[3ex] &=& \dst
\elz\,\f{\ez}{2\pi}\,c_3\,|S(\zz)|.
\ea
\eeq

Our task is thus reduced to establishing the inequality (\ref{eq11}) with constants $c_1$, $c_2$ such that $\sqcc\leq 2\al^{-1}$ (and $\leq\sqrt{\pi^2+4}$ if $\al<1$).

\medskip
{\bf b}. Using the change of variable $\tau\mapsto -\tau$ on $[-u,0]$, we get
\beq{eq13}
 \iuu \f{\tap{} (s+i\tau)}{\tsq}\,d\tau
= 2 e^{\f{i\pi\al}{2}}\;
 \iu \f{(a\tau^{\al+1}+sb\tap{})(\tsq)}{\tsqq}\,d\tau
\eeq
where 
\beq{eq14}
 a=\sin\f{\pi\al}{2},\qquad 
 b=\cos\f{\pi\al}{2}.
\eeq

We will need the integral representations
\beq{eq15}
\f{u^{\al+1}}{1+u^2}\;=\;
\left.\f{\tap{+1}}{\tsq}\right|_{0}^u
\;=\;\iu \f{(\al-1)\tap{+2}+(\al+1)\tap{}}{\tsqq}\,d\tau
\eeq
and
\beq{eq16}
\f{u^{\al+2}}{1+u^2}\;=\;
\iu \f{\al\tau^{\al+3}+(\al+2)\tap{+1}}{\tsqq}\,d\tau.
\eeq

Multiplying the inequality (\ref{eq11}) by $u^\al$ and making the substitutions (\ref{eq13}), (\ref{eq15}), (\ref{eq16}),
we transform (\ref{eq11}) to the equivalent form
\beq{eq17}
\begin{array}{c}
\dst
\left|\iu
\f{a\tau^{\al+3}+s k_{-}\tap{+2}+a\tap{+1}+s k_{+}\tap{}}{\tsqq}\,d\tau
\right|
\\[5ex]
\hspace*{5em} \dst
\leq\;\; \f{1}{2} c_2\,\iu
\f{\al\tau^{\al+3}+(\al+2)\tap{+1}}{\tsqq}\,d\tau,
\end{array}
\eeq
where
$$
k_{\pm}=b-\f{c_1}{2}\,\enpa \,(\al\pm 1).
$$
In order for (\ref{eq17}) to hold for small positive $u$, the coefficient of $\tau^\al$ must equal $0$, so we set
\beq{eq18}
c_1=\f{2b}{\al+1}\,\epa. % \,(\al+1)^{-1}.
\eeq
With this value of $c_1$, the numerator in the l.h.s.\ of
(\ref{eq17}) is real. Clearing the absolute value notation, we rewrite (\ref{eq17}) as a system of two inequalities, in which we leave the least favorable sign of $s$ (so as to make the coefficient of $\tau^{\al+2}$ negative):
$$
\iu \f{(c_2\al\mp 2a)\tau^{\al+3}-{\dst\f{4|b|}{\al+1}}
\tap{+2}+(c_2\al+2c_2\mp 2a)\tap{+1}}{\tsqq}\,d\tau\;\geq\; 0.
$$
Taking the least favorable sign in front of $2a$, we get  
\beq{eq19}
\iu \f{P_2(\tau)\,\tap{+1}}{\tsqq}\,d\tau
\geq 0,
\eeq
where
$$
 P_2(\tau)=\big(c_2\al-2|a|\big)\tau^2
-\f{4|b|}{\al+1}\tau
+\big(c_2\al-2|a|+2c_2\big).
$$
The inequality (\ref{eq19}) with $c_1$ defined by (\ref{eq18}) implies (\ref{eq11}). The rest of the proof amounts to finding an appropriate value of $c_2$. 

\medskip
{\bf c}. It suffices to ensure that the quadratic polynomial $P_2(\tau)$
is nonnegative. 
%
\iffalse

\bigskip
\underline{Comment - March 2011}
Here is an alternative way to proceed. Introduce notation
$$
 F_k(t)=\f{t^k}{(1+t^2)^2}.
$$
Denote, omitting an explicit dependence on $\al$,
$$
G_1(t)=2|a| F_{\al+3}(t)+4\f{|b|}{\al+1} F_{\al+2}(t) +2|a| F_{\al+1}(t)
$$
and
$$
G_2(t)=\al F_{\al+3}(t)+(\al+2) F_{\al+1}(t).
$$
Then
$$
 c_2=\sup_{u>0}\f{\int_0^u G_1(t)\,dt}{\int_0^u G_2(t)\,dt}.
$$
If $u_*$ is the point of maximum, then, by logarithmic differentiation,
$$
  c_2=\f{G_1(u_*)}{G_2(u_*)}.
$$
Therefore,
$$
 c_2\leq \sup_{u>0}\f{G_1(u)}{G_2(u)}. 
$$
Maximization leads to the same quadratic equation as the method in the paper, and to the estimate (\ref{eq20}) for $c_2$.

(End of comment)

\bigskip
\hrule

\bigskip
\fi
%
Let us choose the value of $c_2$ that makes its discriminant equal to zero:
\beq{eq20}
c_2=2\;
%\f{
%|a|(\al+1)+
%\sqrt{a^2(\al+1)^2+\al(\al+2)(b^2(\al+1)^{-2}-a^2)}}{
%\al(\al+2)}
%
\f{|a|(\al+1)^2+\sqrt{a^2+\al(\al+2)}}{\al(\al+1)(\al+2)}.
\eeq
In view of (\ref{eq14}) we have $|a|,|b|\leq 1$, $a^2+b^2=1$,
and the following estimates readily follow: 
\beq{eq21}
\f{|a|}{\al}\;\leq\; \f{c_2}{2}\;\leq\;
\f{|a|(\al+1)+1}{\al(\al+2)}. 
\eeq
The left estimate shows that the coefficient of $\tau^2$ in $P_2(\tau)$ is nonnegative. By our choice of $c_2$, this
leads to (\ref{eq19}) and hence to (\ref{eq11}).

Using the right inequality in (\ref{eq21}) together with
(\ref{eq14}) and (\ref{eq18}), we find $|c_1|^2+c_2^2\leq 4\al^{-2}$, as required.
 
\medskip
{\bf d}. To finish the proof, let us show that if $\al<1$
then one may use the value $c_2=2|a|\al^{-1}$ instead of (\ref{eq20}). With this new choice of $c_2$, the constant $c_3$ in the r.h.s.\ of (\ref{eq12}) becomes
$$
c_3=2\sqrt{\left(\f{\cos(\pi\al/2)}{\al+1}\right)^2+\left(\f{\sin(\pi\al/2)}{\al}\right)^2}
\leq \sqrt{4+\pi^2},
$$
as the theorem claims.
%
\iffalse

\bigskip
\underline{Comment - March 2011}
The right-hand side can be replaced by $sqrt{4+\pi^2(1-\al^2)}$ or even 
$sqrt{4+\pi^2(1-\al^{17/16})}$; the latter exponent is 
fairly close to minimal possible (found numerically).

(End of comment)
 
\bigskip
\hrule

\bigskip

\fi
%
We have to verify the inequality
$$
4\f{|a|}{\al}\,
\iu\f{\tap{+1}}{\tsqq}\,d\tau
-
4\f{|b|}{\al+1}\,
\iu\f{\tap{+2}}{\tsqq}\,d\tau
\geq 0
$$
in the interval $0<\al<1$.
The left-hand side, as a function of $u$, is positive for small $u$ and has a unique critical point (maximum) on $\RR_+$. It remains to check that
$$
\f{|a|}{\al}\,\int_0^\infty\f{\tap{+1}}{\tsqq}\,d\tau
\geq
\f{|b|}{\al+1}\,\int_0^\infty\f{\tap{+2}}{\tsqq}\,d\tau.
$$
Using the substitution $\tau^2=t$ we express the integrals in terms of Euler's Beta function and the last inequality takes the form
$$
\f{|a|}{\al}\, B\left(\f{\al+2}{2},\,\f{2-\al}{2}\right)
\geq
\f{|b|}{\al+1}\,
B\left(\f{\al+3}{2},\,\f{1-\al}{2}\right).
$$
The right and left sides are in fact equal: it follows from (\ref{eq14}) and the identities
$$
B\left(\f{\al+2}{2},\,\f{2-\al}{2}\right)=
\f{\pi\al/2}{\sin(\pi\al/2)},
\qquad
B\left(\f{\al+3}{2},\,\f{1-\al}{2}\right)=
\f{\pi(\al+1)/2}{\cos(\pi\al/2)}.
$$
The proof is complete.

%---------------------------------------------------------
\section{Proof of Theorem~\ref{thm2}}

{\bf a}. The left-hand side of the inequality (\ref{eq8}), likewise the l.h.s.\ of the inequality (\ref{eq5}) in Theorem~\ref{thm1}, vanishes if the closed contour
of integration consisting of $\Gamma$ and the segment $[\ze_0,\zbar_0]$ is used. In the right-hand side of (\ref{eq8}), integration over $\Gamma$ can be replaced by integration over $[\zbar_0,\ze_0]$ since
$$
\iC S_q(\ze) (\lz-\ze)^m\,d\ze
=\int_0^\infty d\Nl\,\iC\f{(\lz-\ze)^m}{(\lam-\ze)^q}\,d\ze,
$$
and the residue of the integrand at $\ze=\lam$ equals $0$ 
%due to 
as
$m\leq q-2$. Therefore (\ref{eq8}) is equivalent to the inequality
\beq{eq22}
\left|\int_0^\infty V_{q,q-1}(\lam)\,d\Nl\right|\;\leq\;
\sum_{m=0}^{q-2} C_m\,
\left|\int_0^\infty V_{q,m}(\lam)\,d\Nl\right|,
\eeq
where
\beq{eq23}
V_{q,m}(\lam)=\eta_0^{q-1-m}\,\int_{\zbar_0}^{\ze_0}\f{(\lz-\ze)^m}{(\lam-\ze)^q}\,d\ze,
\qquad m=0,1,\dots,q-1.
\eeq
The substitutions $\ze=\lz+i\eta_0\tau$,
$\mu=(\lam-\lz)\eta_0^{-1}$ bring (23) to the form
\beq{eq24}
V_{q,m}(\lam)=(-i)^{m-1}\,T_{q,m}(\mu),
\eeq
where
\beq{eq25}
T_{q,m}(\mu)=\int_{-1}^1 \f{\tau^m\,d\tau}{(\mu-i\tau)^q}.
\eeq

\medskip
{\bf b}. Let us study properties of the functions $T_{q,m}(\mu)$.

\medskip
$1^\circ$. The function $T_{q,m}(\mu)$ is even if $q-m$ is even, and odd if $q-m$ is odd. This is verified by changing $\mu$ into $-\mu$ in (\ref{eq25}) and $\tau$ into $-\tau$.

\smallskip
$2^\circ$. If $0\leq m\leq q-2$, then we can write
\beq{eq26}
 T_{q,m}(\mu)=\f{P_{q,m}(\mu)}{(\mu^2+1)^{q-1}},
\eeq
where $P_{q,m}(\mu)$ is a polynomial (even or odd depending on the evenness of $q-m$). Indeed, expanding $\tau^m$ in powers of $\mu-i\tau$, integrating the resulting linear combination of the functions $(\mu-i\tau)^{n-q}$ ($n=0,\dots,m$) with respect to $\tau$, and taking the common denominator, we obtain (\ref{eq26}).

\smallskip
$3^\circ$. From (\ref{eq25}) it is easy to find the asymptotics of $T_{q,m}$ as $\mu\to+\infty$:
\beq{eq27}
\begin{array}{ll}
\dst
T_{q,m}(\mu)=b_{q,m}\,\mu^{-q}\,+\,O(\mu^{-q-2}) &
\;\;\mbox{\rm for $m$ even},
\\[2ex] \dst
T_{q,m}(\mu)=b_{q,m}\,\mu^{-q-1}&
\;\;\mbox{\rm for $m$ odd},
\end{array}
\eeq
where $b_{q,m}\neq 0$.

\smallskip
Comparing to (\ref{eq26}), we see that the polynomial $P_{q,m}$ ($m\leq q-2$) is of exact degree $q-2$ for $m$ even, and of exact degree $q-3$ for $m$ odd.

\medskip
{\bf c}. Let us show that $\{P_{q,m}\}$, $m=0,\dots,q-2$, is a basis in the space of polynomials of degree at most $q-2$. It suffices to verify that the corresponding functions $T_{q,m}$ are linearly independent. Suppose, to the contrary, that some their linear combination is zero:
$$
 L(\mu)=\int_{-1}^1 \f{U(\tau)}{(\mu-i\tau)^q}\,d\tau=0,
$$
where $U(\tau)$ is a polynomial of degree at most $q-2$.
Consider $L$ as an analytic function of complex variable $\mu$. It is regular outside the segment $[-i,i]$, therefore, it equals zero identically. Let $\gamma$ be a closed contour around the segment $[-i,i]$. For all integer $n\geq q-1$ we have
$$
 0=\int_\gamma L(z)\,z^n\,dz
=\int_{-1}^1 U(\tau)\,d\tau\,\int_\gamma\f{z^n\,dz}{(z-i\tau)^q}=\beta_n\int_{-1}^1 U(\tau)\,\tau^{n-q+1}\,d\tau,
$$
and $\beta_n\neq 0$. Hence for any polynomial $\tilde U(\tau)$ 
$$
 \int_{-1}^1 U(\tau)\,\tilde U(\tau)\,d\tau=0.
$$
Taking $\tilde U$ to be the complex-conjugate of $U$ leads to the conclusion $U\equiv 0$, which proves the linear independence of the functions $T_{q,m}$.

\medskip
{\bf d}. Consider first the case of even $q$. The result of ({\bf c}) shows that for some $C_0',\dots,C_{q-2}'$
$$
 \sum_{m=0}^{q-2} C_m' P_{q,m}(\mu)=1+\mu^{q-2},
$$  
or --- cf.\ (26) --- that
\beq{eq28}
 \sum_{m=0}^{q-2} C_m' T_{q,m}(\mu)=\f{1+\mu^{q-2}}{(1+\mu^2)^{q-1}}=:H_{q}(\mu).
\eeq
As $|\mu|\to\infty$, we have
$$
 |T_{q,q-1}(\mu)|\sim |b_{q,q-1}|\,|\mu|^{-q-1}=o\big(H_{q}(\mu)\big).
$$

The function $T_{q,q-1}(\mu)$ is bounded, while $\min H_q(\mu)>0$ on every finite interval. Therefore there exists a positive $C$ such that $|T_{q,q-1}(\mu)|\leq C H_{q}(\mu)$
for all real $\mu\neq 0$. Using (\ref{eq24}) and passing from $T_{q,m}$ to $V_{q,m}$, then integrating with $d\Nl$, we get (\ref{eq22}).

\medskip
{\bf e}. Now consider the case of odd $q$. There exist constants
$C_0'',\dots,C_{q-3}''$ such that
\beq{eq29}
 \sum_{m=0}^{q-3} C_m'' T_{q,m}(\mu)=\f{1+\mu^{q-3}}{(1+\mu^2)^{q-1}}=:H_{q}(\mu).
\eeq
Now $T_{q,q-1}(\mu)$ decays at infinity slower than $H_q(\mu)$. However, as seen from (\ref{eq27}),
$$
 \left|T_{q,q-1}(\mu)-\f{b_{q,q-1}}{b_{q,0}}T_{q,0}(\mu)\right| = O(\mu^{-q-2}) = o(H_q(\mu)),
\quad |\mu|\to\infty.
$$
Here, like previously in ({\bf d}), the left-hand side is a bounded function, while $\min H_q(\mu)>0$ on every finite interval. Therefore there exists a positive $C$ such that for all  real $\mu\neq 0$. 
$$
\left|T_{q,q-1}(\mu)-\f{b_{q,q-1}}{b_{q,0}}T_{q,0}(\mu)\right|\leq C\sum_{m=0}^{q-3}
C_m''\,T_{q,m}(\mu).
$$
Passing from $T_{q,m}$ to $V_{q,m}$, integrating with $d\Nl$, and bringing the term $\mathrm{const}\,|\int_0^\infty V_{q,0}(\lam)\,d\Nl|$ over to the right-hand side, we get (\ref{eq22}).

\smallskip
The proof is finished.

\bigskip
\noindent
{\bf Remark.}\@
Note that due to ({\bf b}$\,1^\circ$) the coefficients of odd functions $T_{q,m}$ in (\ref{eq28}), (\ref{eq29}) are equal to 0. Hence in the right-hand side of (\ref{eq8}) the actual summation is carried over the values of index $m$ for which $q-m$ is even; if $q$ is odd, the value $m=0$ is also included.  

\bigskip
In conclusion I would like to thank M.S.~Agranovich for suggesting the problem and for attention to this work.

\bigskip\bigskip
{\small
Moscow Institute of Electronic Engineering.}

%-----------------------------------------------------
%\section{Appendix}
\subsection*{Remarks added in translation}

%This Appendix includes some remarks not contained in the %published Russian version of the paper.

% Handwritten appendix to master copy of "Pleijel" 
% Dated Oct.1989.

\bigskip\noindent
{\bf 1}. In both Theorems 1 and 2, the function $\Nl$
does not have to be continuous at $\lz$. 
If $\lz$ is a point of discontinuity of $\Nl$, Theorem 1 remains valid without change, while in Theorem 2 the value $N(\lz)$ in the left-hand side can be replaced by any value between $N(\lz-0)$ and $N(\lz+0)$.

\bigskip\noindent
{\bf 2}.
The following theorem stands in the same relation to Theorem 2 as Theorem 1 to Pleijel's inequality (\ref{eq3}).

%(Unlike in Theorem 1, explicit constants are not worked out.) 

\begin{thm} 
\label{thm3}
Let function $N(\lam)$ defined on $\RR_+$ be nondecreasing, equal zero near $\lam=0$, and satisfy condition {\rm (\ref{eq6})}.
For any $\al>0$ and any integer $q=2,3,\dots$
there exist nonnegative constants
$C_0,\dots,C_{q-2}$, depending on $q$ and $\al$, such that
for any $\lam>0$
$$
\begin{array}{l}
\dst
\left| N^{\al}(\lz)\;-\;\f{\al\,B(q,\al)}{2\pi i}\iC
S_q(\ze)\,(\ze-\lz)^{q-1}\,\left(1-\f{\ze}{\lz}\right)^\al\,d\ze
\right|
\\[4ex]
\qquad\qquad\leq\;\;
\dst
\sum_{m=0}^{q-2} C_m\left(\f{\eta_0}{\lz}\right)^\al\,\cdot \eta_0^{q-1-m}\,
\left|\iC S_q(\ze)\,(\ze-\lz)^m\,d\ze\right|,
\end{array}
$$ 
where
$$
B(q,\al)=\f{\Gamma(q)\Gamma(\al)}{\Gamma(q+\al)}.
$$
\end{thm}

\noindent
{\it Proof}: 
Repeat the proof of Theorem 2 replacing the function $T_{q,q-1}(\mu)$ by the function
$$
T_{q,q-1+\al}(\mu)=\int_{-1}^1 \f{\tau^{q-1+\al}}{(\mu-i\tau)^q}\,d\tau.
$$

%---------------------------

\bigskip\noindent
{\bf 3}. An application of Theorem~1 can be found in:
%\bibitem{Agranovich} 
\\[1ex]
Agranovich M.S.
Elliptic operators on closed manifolds. Encycl.\ Math.\ Sci., 63 (Partial Differential Equations VI), Springer-Verlag, 1994, Theorem 6.1.6.
% M.S.A. does not cite applications of the original 
% Malliavin's theorem; only the application of Thm.~1.

\end{document}